\documentclass[12pt]{article}

\usepackage{latexsym}
\usepackage{amsfonts}
\usepackage{amsmath}
\usepackage{amsthm}
\usepackage{amssymb}
\usepackage{amscd}
\usepackage{graphicx}
\usepackage{epstopdf}

\newcounter{num}
\newcommand{\FAT}[1]{\mbox{{$\mathbb{#1}$}}}

\newcommand{\cc}{\FAT{C}}
\newcommand{\CC}{\FAT{C}}

\newcommand{\qq}{\FAT{Q}}

\newtheorem{lem}{\sc Lemma}
\newtheorem{prop}[lem]{\sc Proposition}
\newtheorem{cor}[lem]{\sc Corollary}
\newtheorem{thm}[lem]{\sc Theorem}
\newtheorem{df}[lem]{\sc Definition}

\input{amssym.def}
\input{amssym.tex}
\begin{document}

\title{Graded polynomial identities for   matrices with the transpose involution}
%Graded identities for matrices with transpose involution
\author{Darrell Haile\\
\small Department of Mathematics, Indiana University, Bloomington,
IN 47405 \\
\small haile@indiana.edu\\
\\ %\small and \\
Michael Natapov %
\\
\small Department of Mathematics, Technion-Israel Institute of Technology, Haifa,
Israel 32000 \\
\small natapov@gmail.com}

\date{}                                         %no date appears

\maketitle

%\vspace{1cm}

\begin{abstract}
Let $G$ be a group of order $k$.  We consider the algebra $M_ k(\CC)$ of $k$ by $k$ matrices over the
complex numbers and view it as a crossed product with respect to
$G$   by embedding $G$ in the symmetric group $S_k$ via the
regular representation and embedding $S_k$ in $M_k(\CC)$ in the usual
way. This induces a natural $G$-grading on  $M_k(\CC)$ which we call a
crossed-product grading.  We study the graded $*$-identities for  $M_k(\CC)$
equipped with such a crossed-product grading and the transpose involution.   To each multilinear
monomial in the free graded algebra with involution we associate a directed
labeled graph.  Use of these graphs  allows us to produce a set of generators  for the
$(T,*)$-ideal of identities.  It also leads to new proofs of the results of Kostant and Rowen on the standard identities satisfied by skew matrices.
Finally we determine an
asymptotic formula for the $*$-graded codimension of $M_k(\CC)$.

\end{abstract}

\footnotetext{{\it 2010 Mathematics Subject Classification:} 16S35, 16R10, 16W10, 16W50, 05C30, 05C45.}
\footnotetext{{\it Keywords:} algebra with involution, matrix algebra, graded polynomial identities, graded codimension, directed graphs, Eulerian path.}

\section{Introduction}
In this paper we consider graded $*$--identities on the algebra $M_k(\CC)$ of $k\times k$ matrices over the complex numbers.   A $*$--identity (ignoring the grading)  is a  polynomial over $\qq$ in variables $x_i$ and $x_i^*$ which vanishes on every substitution of matrices in $M_k(\CC)$,  with the condition that if the matrix $A$ is substituted for $x_i$, then $A^*$ (the transpose of A) is substituted for $x_i^*$.  Such identities have been studied by several authors.
If $G$ is a group we say a $\CC$--algebra $B$ is $G$--graded if there is, for each $g\in G$, a  subspace $B_g$ of $B$ (possibly zero)  such that $B=\oplus _{g\in G}B_g$  and for all $g,h\in G$,  $B_gB_h\subseteq B_{gh}$.   Given  a $G$--grading on $M_k(\CC)$ we may consider $G$--graded $*$--identities, that is,  polynomials in weighted variables $x_{i,g}$ and $x_{i,g}^*$  that vanish under all homogeneous substitutions.   In other words we substitute elements from the homogeneous component $M_k(\CC)_g$  for a variable $x_{i,g}$ and if we substitute $A$ for $x_{i,g}$ we must substitute $A^*$ for $x_{i,g}^*$.  We are interested in a particular grading,  the crossed-product grading, on $M_k(\CC)$ given by a group $G$ of order $k$:  To produce this grading we begin by imbedding $G$, via the regular representation, in the group $P_k$ of $k\times k$ permutation matrices.   We then set $M_k(\CC)_g=D_kP_g$  where $D_k$ denotes the set of diagonal matrices and $P_g$ is the permutation matrix corresponding to $g\in G$.  Our main objects of study are  the graded $*$--identities on $M_k(\CC)$ endowed with the $G$--crossed product grading.

Our approach is to use graph theory in a way analogous to our work on  $G$--graded identities on   $M_k(\CC)$,  where the $G$--grading is the crossed-product grading described above.       We start with the free algebra $\qq\{x_{i,g},x^*_{i,g}|i\geq 1, g\in G\}$,  which is a $G$--graded algebra and admits an obvious involution,  also denoted by $*$.   In this free algebra we consider $I(G,*)$,  the ideal of $G$--graded $*$--identities for $M_k(\CC)$.   This ideal is a $(T,*)$--ideal,  which means that it is invariant under  graded endomorphisms and under the involution $*$.    As in the classical case one can show that $I(G,*)$  is generated as a $(T,*)$--ideal by the set of  {\it strongly multilinear}  identities,  where a strongly multilinear  polynomial is a polynomial of the following form: $$\sum_{\pi\in S_n}a_\pi x_{\pi(1),\sigma_{\pi(1)}}^{\epsilon_{\pi,1}} x_{\pi(2),\sigma_{\pi(2)}}^{\epsilon_{\pi, 2}}\cdots x_{\pi(n),\sigma_{\pi(n)}}^{\epsilon_{\pi,n}}$$ where each $
 a_\pi$ is a rational number and for each pair $(\pi,j)$, $\epsilon_{\pi,j}$ is either nothing or  $*$.   The adjective ``strongly"  is to indicate that these monomials are not just multilinear in the variables $x_{i,g}$ and $x_{i,g}^*$  but also in the numerical subscripts, so that we do not allow $x_{i,g}^{\epsilon_i}$ and $x_{i,h}^{\gamma_i}$ to appear in the same monomial unless $g=h$ and $\epsilon_i=\gamma_i$.  As in the case of graded identities without involution we  associate a finite directed graph to each strongly multilinear monomial.   It turns out that two strongly multilinear monomials have the same graph if and only if the difference of the monomials is a graded $*$--identity.

We give three main applications.   The first is a determination of a very simple set of identities that generate  the $(T,*)$--ideal of graded $*$--identities for $M_k(\CC)$ (See Theorem \ref{generators}).
The second application is a new proof of a theorem of Kostant and Rowen on standard identities satisfied by skew matrices.    Our proof uses an analogue of Swan's theorem on paths of length $2k$ in a directed graph with $k$ vertices.    The third application is an asymptotic formula for the codimension of graded $*$--identities for $M_k(\CC)$.  (For a discussion of the notion of codimension and its history see the introduction to our paper \cite{HN}.) To state the result we let $P^G_n$, for each positive integer $n$, be the vector space of strongly multilinear polynomials of degree $n$:

$$P^G_n={\rm span}\{x^{\epsilon_1}_{{\sigma(1)}, g_1}x^{\epsilon_{2}}_{{\sigma(2)}, g_2} \cdots x^{\epsilon_n}_{{\sigma(n)}, g_n}  \  | \ \sigma \in S_n,\  g_1,g_2,\dots , g_n \in G \ \
\hbox{and each   $\epsilon_i$ is $*$ or nothing}\}$$

\noindent This is a vector space of dimension $2^n|G|^nn!$ over  $\qq$.   The $(G,*)$--graded $n$-codimension $c^G_n$ of $M_k(\CC)$ is the dimension of $P^G_n$ modulo the graded $*$--identities:
$$c^G_n = \dim \frac{P^G_n}{P^G_n\cap I(G,*), }.$$

We prove (Theorem \ref{asymptotic}) the following asymptotic result:
$$c^G_n \sim  {k\over 2^{k-1}}k^{2n},$$
a formula interesting both because of its independence of $G$ and its simplicity.  We also determine (Theorem \ref{asymptotic-case-2}) an explicit formula for $c^G_n$ for the case where $|G|$=2.

\section{The graph of a monomial}

In this section we find a set of generators for $I(G,*)$ and show how this set of generators can be understood using graph theory.   We begin by considering the natural group of transformations of the space $P^G_n$.   Let $S_n$ denote the symmetric group and let $C_2$ denote the group of order 2, which we will represent as $\{1,-1\}$ under multiplication.  Let $W_n=S_n \times C_2^n$,  where $C_2^n=C_2\times C_2\cdots C_2$ ($n$ times).   The group $W_n$ acts on $P^G_n$ as follows:  If $\alpha\in W_n$,  we write $\alpha=(\pi,\gamma)$,  where $\pi$ is in $S_n$ and $\gamma=(\gamma_1,\gamma_2,\dots, \gamma_n)$ is in $C_2^n$.  Then if $m=x_{i_1,\sigma_{i_1}}^{\epsilon_{i_1}} x_{i_2,\sigma_{i_2}}^{\epsilon_{i_2}}\cdots x_{i_n,\sigma_{i_n}}^{\epsilon_{i_n}}\in P^G_n$, where $\{i_1,i_2,\dots, i_n\}=\{1,2,\dots, n\}$,  we define $\alpha(m)= x_{\pi(i_1),\sigma_{\pi(i_1)}}^{\delta_{\pi(i_1)}}x_{\pi(i_2),\sigma_{\pi(i_2)}}^{\delta_{\pi(i_2)}}\cdots x_{\pi(i_n),\sigma_{\pi(i_n)}}^{\delta_{\pi(i_n)}}$ where $\delta_j=\epsilon_j$ if $\gamma_j=1$ and $\delta_j$ is the opposite of $\epsilon_j$ (that is $\delta_j$ is nothing if $\epsilon_j=*$ and $\delta_j=*$ is $\epsilon_j$ is nothing) if $\gamma_j=-1$.    So for example if $n=3$ and $\alpha=((1,2,3), (1,-1,-1))$ then $\alpha(x_{2,g_2}^*x_{1,g_1}^*x_{3,g_3})=x_{3,g_3}^*x_{2,g_2}x_{1,g_1}^*$.  We will refer to the elements of $W_n$ as $*$--permutations.

\begin{prop}  The $(T,*)$ ideal $I(G,*)$ is generated by the identities of the following form:

$$f(x_{1,\sigma_1},x_{2,\sigma_2},\dots, x_{n,\sigma_n})=\sum_{\alpha\in W_n}a_\alpha \alpha(x_{1,\sigma_1}x_{2,\sigma_2}\dots, x_{n,\sigma_n}),$$

where $a_\alpha\in \qq$.

\end{prop}

\proof The proof is a standard linearization argument and will be omitted.  \qed

To proceed we need to be more precise about the crossed-product grading on $M_k(\cc)$ obtained from a group $G$ of order $k$.   Let the group $G=\{g_1,g_2,\dots, g_k\}$.  We will parameterize the diagonal elements of an arbitrary diagonal matrix $E\in M_k(\cc)$  using this ordering on $G$,  so the $(i,i)$ entry will be labeled $e_{g_i}$.   Then for each $g\in G$ we can choose a permutation matrix $P_g\in M_k(\cc)$  such that for every diagonal matrix $E$  we have $P_gEP_g^{-1}=E^g$ where the $E^g$ is the diagonal matrix with $(i,i)$ entry $e_{gg_i}$.

Now let $S=\{t_{j, g_i} \ | \ 1\leq j\leq n, 1\leq i\leq k\}$ be a set of commuting indeterminates and let $L=\qq[S]$.    As above we have an action of $G$ on the set of diagonal matrices in $M_k(L)$.  We consider the diagonal matrices $D_1,D_2,\dots, D_n$,  where $D_j$  is the diagonal matrix with $(i,i)$ entry $t_{j,g_i}$.   So the $(i,i)$ entry of the matrix $D_j^g$ is $t_{j,gg_i}$.

\begin{lem}\label{diagonal} Let $(\sigma_{1,1}, \sigma_{1,2},\dots, \sigma_{1,n}),  (\sigma_{2,1}, \sigma_{2,2},\dots, \sigma_{2,n}), \dots,  (\sigma_{m,1}, \sigma_{m,2},\dots, \sigma_{m,n})$ be distinct $n$--tuples of element of $G$.   Then the $m$ diagonal matrices
$$D_1^{\sigma_{1,1}}D_2^{\sigma_{1,2}}\cdots D_n^{\sigma_{1,n}},\ \  D_1^{\sigma_{2,1}}D_2^{\sigma_{2,2}}\cdots D_n^{\sigma_{2,n}},\dots, D_1^{\sigma_{m,1}}D_2^{\sigma_{m,2}}\cdots D_n^{\sigma_{m,n}}$$ are linearly independent over $\qq$.

\end{lem}

\proof  Let $E_i$ be the diagonal matrix $D_1^{\sigma_{i,1}}D_2^{\sigma_{i,2}}\cdots D_n^{\sigma_{i,n}}$.  The $(1,1)$ entry of $E_i$ is \\ $t_{1,\sigma_{i,1}g_1}t_{2,\sigma_{i,2}g_1}\cdots t_{n,\sigma_{i,n}g_1}$.  If $j$ is different from $i$ then the $(1,1)$ of $E_j$ entry is equal to the $(1,1)$ entry of $E_i$ if and only if  $\sigma_{i,r}g_1=\sigma_{j,r}g_1$ for all $r, 1\leq r\leq n$.  This implies that $\sigma_{i,r}=\sigma_{j,r}$ for all $r$, so the two $n$-tuples are not distinct.  Therefore the $(1,1)$ entries of the $m$ diagonal matrices $E_1,E_2,\dots, E_m$ are distinct monomials and hence linearly independent over $\qq$.  It follows that $E_1,E_2,\dots, E_m$ are also linearly independent over $\qq$.\qed

\begin{df}\label{df:pathtransformation}
Let $m=x_{1,\sigma_1}^{\epsilon_1} x_{2,\sigma_2}^{\epsilon_2}\cdots x_{n,\sigma_n}^{\epsilon_n}$   and  let $\alpha=(\pi,\gamma)$   be a $*$--permutation.  Let $\gamma=(\gamma_1,\dots, \gamma_n)$.    We will call $\alpha$ a \underbar{path} \underbar{transformation} of $m$ if the following conditions hold:

(1)  $\sigma_1^{a_1}\sigma_2^{a_2}\cdots \sigma_{n}^{a_{n}}=\sigma_{\pi(1)}^{c_{\pi(1)}}\sigma_{\pi(2)}^{c_{\pi(2)}}\cdots \sigma_{\pi(n)}^{c_{\pi(n)}}$, where  $a_r=-1$ when $\epsilon_r = *$ and 1 otherwise ,  and $c_r=\gamma_r a_r$.

(2) For all $i$, $1\leq i\leq n$,  if $j$ is the unique integer such that $\pi(j)=i$,  then  $\sigma_1^{a_1}\sigma_2^{a_2}\cdots \sigma_{i-1}^{a_{i-1}}\sigma_i^{b_i}=\sigma_{\pi(1)}^{c_{\pi(1)}}\sigma_{\pi(2)}^{c_{\pi(2)}}\cdots \sigma_{\pi(j-1)}^{c_{\pi(j-1)}}\sigma_{\pi(j)}^{d_{\pi(j)}}$, where  $a_r=-1$ when $\epsilon_r = *$ and 1 otherwise,  $b_i=-1$ when $\epsilon_i=*$ and 0 otherwise,  $c_r=\gamma_r a_r$ when $\gamma_r = *$ and 1 otherwise, and  $d_i=-1$ if $\gamma_ia_i=-1$ and 0 otherwise.

\end{df}
\bigskip

\noindent Remark: The reason for the name path transformation will become clear when we introduce the graph of a monomial.   Also we will see later that condition (2) implies condition (1).

\bigskip

We will call two monomials $m,r$ in $\qq\{x_{i,g},x^*_{i,g}|i\geq 1, g\in G\}$  \underbar{equivalent} if $r=\alpha(m)$ for some path transformation $\alpha$ of $m$.   This is clearly an  equivalence relation on the set of monomials in each $P^G_n$.

\begin{prop} Let $f(x_{1,\sigma_1},x_{2,\sigma_2},\dots, x_{n,\sigma_n})=\sum_{\alpha\in W_n}a_\alpha \alpha(x_{1,\sigma_1}x_{2,\sigma_2}\cdots x_{n,\sigma_n}),$ for some $a_\alpha\in \qq$.
Let $f=f_1+f_2+\cdots +f_t$ be the decomposition of $f$ into the sums over equivalent monomials.   The $f$ is a graded $*$--identity for $M_k(\cc)$ if and only if each $f_i$ is a graded $*$--identity  and a given $f_i$ is a graded $*$--identity if and only if the sum of its coefficients is zero.
\end{prop}

\proof  In any evaluation of $f$ we must substitute a matrix of the form $DP_\sigma$ for a variable  $x_\sigma$ of group weight $\sigma$.
Note that when you make the  substitution $DP_\sigma$ for a term of the form $x_\sigma^*$ you get  $(DP_\sigma)^*=P_\sigma^*D^*$.  But the transpose of a permutation matrix is the inverse of that matrix, so $P_\sigma^*=P_{\sigma^{-1}}$, and diagonal matrices are symmetric, so  $D^*=D$.   Therefore  $(DP_\sigma)^*=P_{\sigma^{-1}}D=D^{\sigma^{-1}}P_{\sigma^{-1}}$.   Now suppose $f$ is an identity and substitute $D_iP_{\sigma_i}$ for $x_{i,\sigma_i}$.   Because we are working over $\cc$,  we may assume that the $nk$ entries in the diagonal matrices $D_1, D_2,\dots, D_n$ are algebraically independent over $\qq$.     Now let $m=x_{1,\sigma_1}^{\epsilon_1} x_{2,\sigma_2}^{\epsilon_2}\cdots x_{n,\sigma_n}^{\epsilon_n}$  be a monomial appearing in $f$  and choose some  $i$, $1\leq i\leq n$.  (To be absolutely general we would need to use $x_{j_i}$ instead of $x_i$ but the argument is the same in that case and the notation less unwieldy in this case).  The contribution to the diagonal part  of the evaluation of the monomial $m$ that comes from substituting $D_iP_{\sigma_i}$  for $x_{i,\sigma_i}$ is  $$D_i^{\sigma_1^{a_1}\sigma_2^{a_2}\cdots \sigma_{i-1}^{a_{i-1}}\sigma_i^{b_i}}$$ where $a_r=-1$ when $\epsilon_r = *$ and 1 otherwise, $b_i=-1$ when $\epsilon_i=*$ and 0 otherwise.
Now let $j=\pi^{-1}(i)$.  Any other monomial that appears in $f$ is of the form  $\alpha(m)$ for some $\alpha \in W_n$.  If $\alpha=(\pi,\gamma)$ and $\gamma=(\gamma_1,\dots, \gamma_n)$,   then the contribution to the diagonal part of $\alpha(m)$ that comes from substituting $D_iP_{\sigma_i}$ is $$D_i^{\sigma_{\pi(1)}^{c_{\pi(1)}}\sigma_{\pi(2)}^{c_{\pi(2)}}\cdots \sigma_{\pi(j-1)}^{c_{\pi(j-1)}}\sigma_{\pi(j)}^{d_{\pi(j)}}}$$ where $c_r=\gamma_ra_r$ when $\gamma_r = *$, and  $d_i=-1$ when $\gamma_i=-1$ and when $\gamma=1$. Because we are selecting our diagonal matrices with algebraically independent entries, we are in the situation of Lemma \ref{diagonal}.   We therefore see that   if this evaluation is zero it follows that the sum of the monomials in which, for all $i$,  the group element applied to $D_i$  is the same must be an identity.   Because we are looking for  $G$--graded identities we must also have that the sum of the  monomials that give the same  group weight must be an identity.   But these two conditions are precisely the conditions of equivalence of monomials:  Condition (1) is the condition to have the same group weight and condition (2) is the condition that the group element applied to $D_i$ is the same for all $i$. So we have that if $f$ is an identity then each $f_i$ vanishes upon substitution with algebraically independent entries for the diagonal elements.    By specialization it follows that in fact each $f_i$ is itself an identity.  But for a given $f_i$  the substitution described above has the same diagonal part and permutation matrix part for each monomial.   Hence $f_i$ is an identity if and only if the sum of the coefficients is zero.\qed

\begin{cor}\label{cor:binomials}
The space  $I(G,*)\cap P^G_n$ is spanned by the elements  $m - \alpha(m)$ where $m$ is any monomial in  $P^G_n$  and $\alpha$ is a path transformation for $m$.
\end{cor}

\proof Clear \qed

We now introduce the graph of a strongly multilinear monomial.   As in the case of graded identities without an involution,  the kind of graphs we consider will be finite directed graphs, with labels on the vertices and the edges.  Every edge has a direction.  There may be several edges in both directions between two given vertices and there may be edges with the same beginning and ending vertex. We begin with a finite group $G$ of order $k$.   For each strongly multilinear monomial
$m=x_{i_1,\sigma_1}^{\epsilon_1} x_{i_2,\sigma_2}^{\epsilon_2}\cdots x_{i_n,\sigma_n}^{\epsilon_n}$ in the free algebra $\qq\{x_{i,g},x^*_{i,g}|i\geq 1, g\in G\}$  we first form a preliminary graph:   There are $k$ vertices labeled by all of the elements of the group.   If there are no starred variables then the preliminary graph is just the graph of the monomial we considered in the case of grading without involution (and this will be the final graph in this case, too).    We remind the reader of that construction:   There is an edge labeled $i_1$ from the vertex labeled $e$ to the vertex $\sigma_1$ and, for $j>1$,  an edge labeled $i_j$ from the vertex labeled $\sigma_1\sigma_2\cdots \sigma_{j-1}$ to the vertex  $\sigma_1\sigma_2\cdots \sigma_{j-1}\sigma_j$ .  In other words the graph is really a directed path through the vertices starting at $e$ and passing successively through $\sigma_1,\sigma_1\sigma_2,\sigma_1\sigma_2\sigma_3$, and so on, ending at $\sigma_1\sigma_2\cdots \sigma_n$.   This is in fact a directed Eulerian path through the graph (that is, it uses each edge exactly once).  Now, if there are starred variables, then the preliminary graph is obtained by starting with the   basic graph of the monomial obtained  by changing each starred variable $x_\sigma^*$ to $x_{\sigma^{-1}}$ and adding a star to each edge that came from a starred variable.   Then the final graph is obtained by reversing each starred edge and removing the star.   Note that when we reverse the edge we are changing the weight of that edge from $\sigma^{-1}$ back to  $\sigma$.  This graph also has an Eulerian path: The edges are numbered $1,2,3\dots $,  but this path is not necessarily directed; some of the edges may be reversed.   Moreover one can read off the monomial from this path:  The order of the numbering of the variables is given by the ordering of the edges in the path, the weights are the weights determined by those edges, and the ``exponent" is $*$ if the path reverses the direction of the edge and nothing if  it goes in the same direction.

\bigskip
\noindent Example:  Let $G=\langle \sigma \rangle$ be the cyclic group of order 3.   Consider the monomial $x_{1,\sigma}^*x_{2,\sigma}x_{3,\sigma}x_{4,e}x_{5,\sigma^2}^*$.   Here are the (a) preliminary  and (b) final graphs:

\graphicspath{{./}}
\begin{figure}[h]
\caption {}\label{fig:Stardef}
$$\includegraphics[width=120mm]{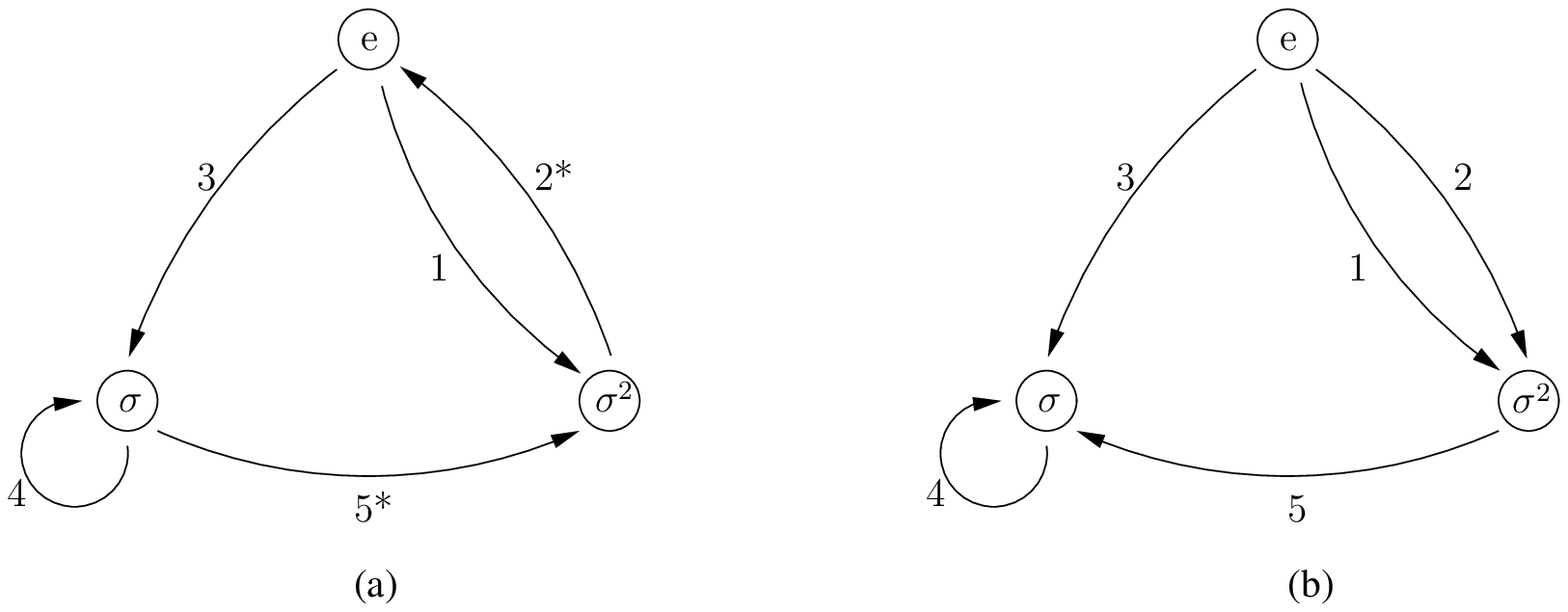}$$
%  \caption{Let $G=\{e,\sigma, \sigma^2\}$ be the cyclic group of order 3, generated by $\sigma$. Given are (a) the preliminary graph, and (b) the final graph of the monomial $x_{1,\sigma}^*x_{2,\sigma}x_{3,\sigma}x_{4,e}x_{5,\sigma^2}^*$.}\label{fig:Stardef}
\end{figure}

\begin{thm} \label{thm:the same graph}
Let $m$ be a strongly multilinear monomial and let $\alpha(m)$ be a $*$--permutation of $m$.   Then $\alpha$ is a path transformation of $m$ if and only if $m$ and $\alpha(m)$ have the same graph.
\end{thm}

\proof
We will in fact prove that $m$ and $\alpha(m)$ have the same graph if and only if  condition  (2) of the definition of path transformation is satisfied.   It follows from the discussion following the definition of the graph associated to a monomial that to say that the monomials $m$ and $\alpha(m)$ have the same graph is the same as saying that in the graph of $m$ there is an Eulerian path starting at $e$  such that the resulting monomial (where you use $x_{\sigma_i}^*$   if your path reverses the direction of the $i$-th edge) is the monomial $\alpha(m)$.    By the definition of the graph,  in any path from $e$ to a vertex $\sigma$  the product of the weights along the path (where you replace $g$ by $g^{-1}$ if you reverse the direction of an edge)  is equal to $\sigma$.
So assume $m$ and $\alpha(m)$ give the same graph,  so that $\alpha(m)$ determines a path on the graph of $m$.  We \underbar{claim} that for each $i$ the equality of condition (2)  arises from two paths from $e$ to the beginning point of the $i$-th edge.    To prove the claim we consider the four cases arising for the possible choices for $b_i$ and $d_i$.   If $b_i=d_i=0$, then because $\pi(j)=i$ we see that the claim is clear:  both sides represent the products of weights in a path from $e$ to the beginning point of the $i$-th edge.  If $b_i=0$ and $d_i=-1$  then the left hand side is again the product of the weights of the edges ending at the beginning point of the $i$-th edge.  On the right hand side because the last term is $\sigma_{\pi(j)}^{-1}$ ,  the path reverses the direction of $\sigma_i$ (recall $\pi(j)=i$ and so the path ends at the \underbar{beginning} point of the  $i$-th edge, as claimed).   The other two cases are similar.    The converse works the same way:   We need to show that the  edge  labeled $i$ in the graph of $m$ is in the same place in the graph of $\alpha(m)$ (by definition of the graph of a monomial,  the edge labeled $i$ has weight $\sigma_i$ in both graphs).   It is enough to show that the edge labeled $i$ has the same beginning vertex in both graphs.    But for each $i$  the equation of (2)  says exactly that:   The product on the left is the group value of the beginning vertex of the edge labeled $i$ in the graph for $m$ and the product on the right is the group value of the beginning vertex of the edge labeled $i$ in the graph of $\alpha(m)$. \qed

\begin{cor}  In the definition of path transformation,  Definition \ref{pathtransformation},  condition (2) implies condition (1).
\end{cor}

\proof  In the proof of the theorem we proved that $m$  and $\alpha(m)$ have the same graph if and only if condition (2) is satisfied.    But if they have the same graph then the ending vertices of the path for $m$ and the path for $\alpha(m)$  must be the same:  it is either $e$ or the only vertex other that $e$ at which there are an odd number of edges either beginning or ending there.  The group value of this final vertex is then the common value for condition (1).\qed

Here is another way to  state the theorem:  given the graph of $m=x_{i_1,\sigma_1}^{\epsilon_1} x_{i_2,\sigma_2}^{\epsilon_2}\cdots x_{i_n,\sigma_n}^{\epsilon_n}$, one obtains all the monomials $r$ such that $m-r$ is a graded $*$--identity as follows:  Look for (nondirected) Eulerian paths through the graph of $m$ starting at e.   Given such a path you associate the monomial  $r$ which has the same variables as $m$ in the order determined by the new path and for which $x_{i_j,\sigma_j}$ now has the same exponent $\epsilon_j$  as $m$ if the direction of the corresponding edge is the same as for $m$ and has the other choice for $\epsilon_j$ (that is, if the old $\epsilon_j$ was nothing the new one is * and if the old was * the new is nothing) if the direction of the edge is reversed.

\medskip

Here is an example:

\medskip

On 2 by 2 matrices with the $C_2$ crossed-product grading you start with the graph of the monomial $x_{1,\sigma}x_{2,\sigma}x_{3,\sigma}$.  So your path is edge 1 followed by edge 2 followed by edge 3.  You look at a new path (ignoring directions of the 3 edges) going edge 1 followed by edge 3 followed by edge 2.  The new monomial is $x_{1,\sigma}x_{3,\sigma}^*x_{2,\sigma}^*$ because you reversed the directions of edges 2 and 3.  The difference $x_{1,\sigma}x_{2,\sigma}x_{3,\sigma}-x_{1,\sigma}x_{3,\sigma}^*x_{2,\sigma}^*$
is a graded identity.

\medskip

Here  is another example:

\medskip

On 6 by 6 matrices with the $S_3$ crossed-product grading you start with the monomial $x_{1,\sigma}x_{2,\sigma\tau}x_{3,\sigma\tau^2}x_{4,\tau^2}x_{5,\sigma\tau}x_{6,\sigma\tau^2}$.  If you draw the graph of this monomial you see that the endpoint is $\sigma\tau$.  There are many other paths from $e$ to $\sigma\tau$ (ignoring the directions of the edges).  For example there is a path $1,2,5^*,4^*,6^*,3$ where the $*$ means that the direction for that edge is reversed.  So we get the binomial graded identity  $x_{1,\sigma}x_{2,\sigma\tau}x_{3,\sigma\tau^2}x_{4,\tau^2}x_{5,\sigma\tau}x_{6,\sigma\tau^2} - x_{1,\sigma}x_{2,\sigma\tau}x_{5,\sigma\tau}^*x_{4,\tau^2}^*x_{6,\sigma\tau^2}^*x_{3,\sigma\tau^2}$.

\bigskip
\section{Generators for the ideal of identities}

\medskip
\begin{thm} \label{generators} Let $G$ be a group of order $k$. The ideal of graded $*$-identities of $M_k(\cc)$ endowed with the $G$--crossed-product grading is generated as a  $(T,*)$--ideal by the following elements:

%\medskip

(1) $x_{i,e}x_{j,e}-x_{j,e}x_{i,e}$ for all $i,j\geq 1$.

(2) $x_{i,e}-x_{i,e}^*$ for all $i\geq 1$

\end{thm}

Remarks:  (1) Notice that these generators are independent of the group; they refer to the $e$--component only.

(2) Bahturin and  Drensky have shown that the $T$--ideal of $G$--graded identities of $M_k(\cc)$  (no involution) is generated by the identities of type (1)  and identities of the following type:  $x_{1,\sigma}x_{2,\sigma^{-1}}x_{3,\sigma}-x_{3,\sigma}x_{2,\sigma^{-1}}x_{1,\sigma}$, for all $\sigma\in G$.   Because every $G$--graded identity without involution is also an identity with involution,  if our result is to be true we must be able to generate this identity from our two kinds.   In fact we can:
Because $x_{2,\sigma^{-1}}x_{3,\sigma}$ has weight $e$, we can apply the second generator to get that $x_{1,\sigma}x_{2,\sigma^{-1}}x_{3,\sigma}-x_{1,\sigma}x_{3,\sigma}^*x_{2,\sigma^{-1}}^*$ is an identity.  Then because $x_{1,\sigma}x_{3,\sigma}^*$ has weight $e$ we get that $x_{1,\sigma}x_{3,\sigma}^*x_{2,\sigma^{-1}}^*-x_{3,\sigma}x_{1,\sigma}^*x_{2,\sigma^{-1}}^*$ is an identity.  Finally because $x_{1,\sigma}^*x_{2,\sigma^{-1}}^*$  has weight $e$ we get that $x_{3,\sigma}x_{1,\sigma}^*x_{2,\sigma^{-1}}^*-x_{3,\sigma}x_{2,\sigma^{-1}}x_{1,\sigma}$ is an identity.  Adding these three gives that $x_{1,\sigma}x_{2,\sigma^{-1}}x_{3,\sigma}-x_{3,\sigma}x_{2,\sigma^{-1}}x_{1,\sigma}$ is an identity, as desired.

\proof  As we have seen in Corollary \ref{cor:binomials},  the graded $*$--identities are generated by the binomial identities $m-\alpha(m)$, where $m$ is any strongly multilinear monomial and $\alpha$ is a path transformation for $m$.   By the remarks above, given $m$, the $\alpha(m)$ where $\alpha$ is a path transformation for $m$ correspond to the  paths in the graph of $m$.   The two kinds of generators can be viewed as basic ways of changing from one path to another in a given graph.   The type (1) identities say that given a path with two consecutive path segments each of which forms a loop at some given vertex $\sigma$ then one can switch the order of those two path segments (clearly that gives a new path through the graph).   The type (2) identities say that given a path with a path segment which is a loop at some vertex $\sigma$ one can form a new path by reversing all the edges of that loop.  What we need to prove then is that one can obtain every  path through the vertices of the graph for $m$ by starting with the path corresponding to $m$ and applying these two kinds of basic moves.

To prove this we begin with a special case, which is the graphical version of a generalization of remark (2) above.   We begin with a path which contains vertices $P$ and $Q$  for which there are  two different path segments $\rho_1$ and $\rho_2$ going from $P$ to $Q$.   Clearly if we switch these two path segments we obtain another path through the vertices.  We \underbar{claim} that this new path can be obtained from the original path by basic moves:   Assuming $\rho_1$ precedes $\rho_2$ in the original path we let $\gamma$ denote the segment of the path between $\rho_1$  and $\rho_2$, so we have a segment of our original path that is $\rho_1\gamma\rho_2$.  The path segment $\rho_1\gamma$ is a loop from $P$ to itself and so we can reverse that loop (applying a basic move of type (2)) and get a new path with the segment  $\gamma^*\rho_1^*\rho_2$, where the star indicates that we are reversing the direction of that path segment.   Now in this new path the  segment $\rho_1^*\rho_2$ is also a loop from $P$ to itself,  we can reverse that loop to get a new path with the segment $\gamma^*\rho_2^*\rho_1$.   Finally in this path the segment $\gamma^*\rho_2^*$ is a loop from $Q$ to itself so we can reverse that loop to get a new path with the segment $\rho_2\gamma\rho_1$,   so we have switched $\rho_1$ and  $\rho_2$, which is what we wanted to accomplish.  We will refer to this kind of move as a \underbar{switch} move.

The rest of the proof is quite similar to the proof of the analogous fact in Haile and Natapov \cite[Corollary 6]{HN}). Again we start with two paths through the graph for $m$,  the original path that determines $m$ , and the new path.    Without loss of generality we may assume that the graph of $m$  has no reversed edges,  that is, the monomial $m$ has no starred variables.   We will list its edges in order as edge 1, edge 2, and so on.    In the new path  edge 1 may be in a different position and may be reversed (which we will denote edge $1^*$).  We will prove that by applying basic moves (or  the switch move which we have shown comes from basic moves)  we can move edge 1  (or edge $1^*$)   back to being the first edge in the path and not being reversed.   It will follow by induction that we can move the new path back to the original path by basic moves.

 First assume that in the new path  edge 1 is reversed, so that edge $1^*$ appears.    There is  a initial loop from $e$ to itself, the last edge of which is  edge $1^*$.   Reversing this loop puts edge 1 at the beginning of the path, so we are done in that case.   From now on we can assume the new path contains edge 1.

Now assume edge 2 or edge $2^*$  appears before edge 1 in the path.   There is then a loop from e to itself that ends with the edge appearing before edge 1 in this path.    If edge $2^*$ appears we can reverse this loop to get  edge $2$  appearing before edge 1  instead and so we have reduced ourself to the case where edge 2 appears before edge 1.   In that case there are two path segments from e to the endpoint of edge 1.   The first is the initial path segment that ends just before edge 2 and the second is edge 1 itself.   If we apply the switch move to these two path segments then we put edge 1 at the beginning of the path, as desired.

So we may assume edge 2 appears after edge 1.    If edge 3 or edge $3^*$  appears before edge 1,   then as above we can first reduce to the case where  edge 3 appears before edge 1.  We then have two path segments from e to the beginning vertex of edge 3.   The first is the initial path segment that ends just before edge 3 and the second is the path segment that begins at the beginning of edge 1 (that is, at e) and ends just after edge 2.    Switching these two segments puts edge 1 back at the beginning of the path.    Continuing in this way,  whenever we have edge k or edge k$^*$ appearing before edge 1 and edge 1 appearing before edge k-1 we can apply basic moves to  put edge 1 at the beginning of the path.   This finishes the proof.\qed

\bigskip

\medskip
\section{Theorem of Kostant and Rowen}

\medskip

We begin by defining the sign of a $*$--permutation:    For each $n$,  we have the product homomorphism  $\mu$   from $C_2^n$ to $C_2$ given by $\mu((\gamma_1,\gamma_2,\dots,  \gamma_n)=\gamma_1\gamma_2\cdots \gamma_n$.   We then define a homomorphism $sign$  from $W_n=S_n\times C_2^n$ to $C_2$ by $sign(\pi,\gamma)=sgn(\pi) \mu(\gamma)$,  where $sgn(\pi)$ is the usual sign homomorphism on $S_n$.   So for example if $m=x_{1,\sigma_1}x_{2,\sigma_{2}}^*x_{3,\sigma_3}$
and $\alpha(m)=x_{2,\sigma_2}x_{1,\sigma_1}^*x_{3,\sigma_3}^*$ then the permutation is odd and the $\gamma=(-1,-1,-1)$, so the sign of $\alpha$ is $(-1)\cdot(-1)^3=1$.

In terms of the graph of $m$  we have seen that each path transformation $\alpha$ corresponds to a path from $e$ to the endpoint of the path for $m$ with possible reversed edges.  The sign of $\alpha$ is then the sign of the permutation of the edges multiplied by $(-1)^k$ where $k$ is the number of reversed edges in the path corresponding to $\alpha$.  We call $\alpha$ \underbar{even} if its sign is 1 and \underbar{odd} if its sign is -1.

Now note that if $f(Y_1,Y_2,\dots, Y_n)$ is a (nongraded) identity  on the skew elements  ($x=-x^*$) of $M_k(\CC)$ then the polynomial $f(X_1-X_1^*, X_2-X_2^*,\dots, X_n-X_n^*)$ is a $*$--identity of on all of $M_k(\CC)$ and this new polynomial is multilinear if $f$ is multilinear.  In particular if you start with a standard polynomial $s(Y_1,Y_2,\dots, Y_n)(=\sum_{\pi\in S_n}sign(\pi)Y_{\pi(1)}Y_{\pi(2)}\cdots Y_{\pi(n)})$ which is an identity on skew elements (for example if $n$ is greater than the dimension of the space of skew elements), then $S(X_1,X_1^*, X_2, X_2^*, \dots, X_n,X_n^*)=s(X_1-X_1^*, X_2-X_2^*,\dots, X_n-X_n^*)$ is an identity on $M_k(\CC)$.  It follows that for every choice of $\sigma_1,\sigma_2,\dots, \sigma_n\in G$ the polynomial $S(X_{1, \sigma_1},X_{1, \sigma_1}^*, X_{2, \sigma_2},X_{2, \sigma_2}^*,\dots, X_{n, \sigma_n},X_{n, \sigma_n}^*,)$ is a graded identity.  If we recall that two monomials are said to be equivalent if one is a path transformation
of the other and we divide the monomials appearing in $S(X_{1, \sigma_1},X_{1, \sigma_1}^*, X_{2, \sigma_2}, X_{2, \sigma_2}^*, \dots, X_{n, \sigma_n},X_{n, \sigma_n}^*,)$ into equivalence classes then the sum of the coefficients of the monomials in each equivalence class must be zero.  But if you choose a monomial $m$ in a class then the coefficient of $\alpha(m)$ is just the sign of $\alpha$ times the coefficient of $m$, so we obtain the following statement about the graph of $m$:  The number of paths starting at e and ending at the endpoint of $m$ must be even and half come from even path transformations of $m$ and half from odd path transformations.  We conclude that there must be a theorem about labeled, directed graphs of the following form:
For every labeled, directed graph on $k$ vertices with a path of "sufficiently large" length n there are an even number of path transformations of this path, half odd and half even.   By what we have just seen, because  the dimension of the space of skew symmetric matrices in $M_k(\CC)$  is $(1/2)k(k-1)$), the number $n= 1+(1/2)k(k-1)$ is sufficiently large.  But in fact a much better result is known.   It was shown by Kostant \cite{Ko} that if $k$ is even, the space of $k\times k$ skew matrices satisfies the standard identity $s_{2k-2}$ and Rowen \cite{Ro} proved that the same result holds for $k$ odd.  We therefore have the following purely graph-theoretic result.

\begin{thm}\label{swanlike} Let $\Gamma$ denote a directed graph on $k$ vertices containing a (not necessarily directed)  path $\rho$ of length $n$ between the vertices $P$ and $Q$.   If $n$ is at least $2k-2$ then there are an even number of paths from $P$ to $Q$ (using the same set of $n$ edges),  half even transformations of $\rho$ and half odd.
\end{thm}

\begin{figure}[h]
 \caption{}\label{fig:Counterexample}
$$\includegraphics[width=100mm]{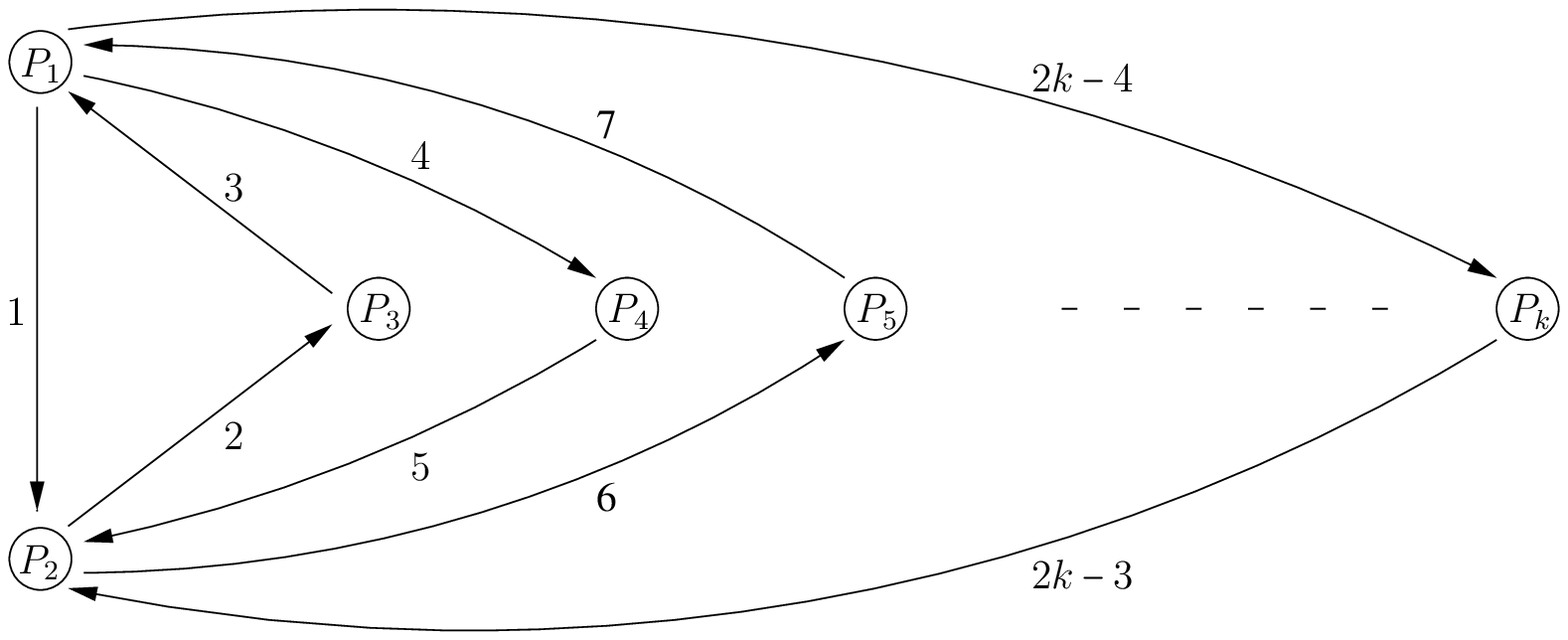}$$
 % \caption{}\label{fig:Counterexample}
\end{figure}

Of course this is analogous to Swan's result \cite{Sw, Sw2} on directed paths. Kostant and Rowen also show that their result is best possible.    We give an example to show the same for the graph result:   In the following graph (Figure \ref{fig:Counterexample}) there are $k$ vertices and a path through all of the $2k-3$ edges, but every path transformation is even.

It is natural to seek a graph theoretic proof of this result.  Of course such a proof will then give a new proof of the theorem of Kostant and Rowen.    The rest of this section is devoted to such a proof.

The idea of the proof is to follow Swan's strategy:   Reduce to the case of balanced graphs and then use induction.   Complications arise because of the possibility of reversed edges.  Recall that in a graph $\Gamma$ the \underbar{order} of a vertex $P$ is the total number of edges coming into or leaving $P$ (an edge that forms a loop at $P$ is counted twice).

\begin{prop}\label{two vertices}
Let $\Gamma$ be a balanced graph with $k$ vertices and $2k-2$ edges.   There are at least two vertices $P$ with order two.

\end{prop}
\proof  By changing the direction of edges we change to a new graph $\Gamma^\prime$ with has a directed loop through the vertices.   For each vertex $A$ of $\Gamma^\prime$  let $n_A$ be the number of edges with endpoint at $A$.   Then clearly $2k-2=\sum_An_A$.  Let $r$ be the number of vertices $A$ such that $n_A=1$.   Then we have the inequality $2k-2\geq r+2(k-r)$.  It follows that $r\geq 2$.    It then follows that at these two vertices in the graph $\Gamma$ the total number of edges going in out of each vertex is two. \qed

Let $\Gamma$ be a balanced graph with $k$ vertices and $2k-2$ edges.  Let $\alpha$ be an Eulerian path for $\Gamma$.   If $P$ is a vertex with order two we will say an Eulerian path $\beta$ has the same direction as $\alpha$ at $P$ if the $\beta$ does not reverse the direction of $\alpha$ at the edges going in or coming out of $P$.

\begin{prop}\label{loops} Let $\Gamma$ be a balanced graph with $k$ vertices and $2k-2$ edges.  Let $\alpha$ be a Eulerian path for $\Gamma$.   If $P$ is a vertex of order two,  there are an even number of Eulerian paths $\beta$ that  have the same direction as $\alpha$ at $P$,  half of them even transformations of $\alpha$ and half of them odd.
\end{prop}
\proof
We may assume that $\alpha$ is a directed Eulerian path for $\Gamma$.   In that case $n_P=1$.  The proof is by induction on $k$.   By Proposition \ref{two vertices} there is a second vertex $Q$ with $n_Q=1$.   We consider several cases:

\begin{figure}[h]
 \caption{}\label{fig:swan1}
$$\includegraphics[width=90mm]{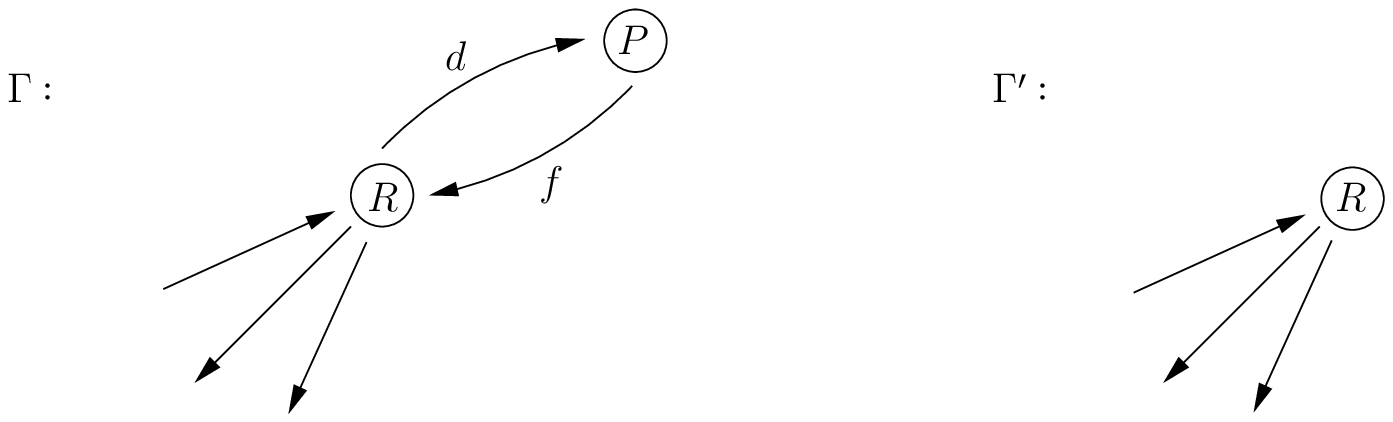}$$

\end{figure}

Case 1)   (See Figure \ref{fig:swan1})  The initial vertex $R$ for the edge $d$ ending at $P$ is the same as the terminal vertex for the edge $f$  starting at $P$.  In this case we look at the graph $\Gamma^\prime$ and the directed path $\alpha^\prime$  obtained by removing $P$ and the two edges $d$ and $f$.   If $R=Q$, then $k=2$ and it is easy to check that the result is true.   So we may assume $R\not= Q$.   In that case the induction hypothesis applies to $\Gamma^\prime$  and so there are an even number of paths for $\Gamma^\prime$  having the same direction at $Q$ as $\alpha^\prime$,  half odd and half even path transformations of  $\alpha^\prime$.   It follows that the total number of paths on $\Gamma^\prime$ is also even,  half odd and half even,   because the others are obtained by simply reversing the direction of those that have the same direction at $Q$ as $\alpha^\prime$.    But these paths are clearly in one to one correspondence with the paths on $\Gamma$ that h
 ave
the same direction at $P$ as $\alpha$,  so we are done in this case.

\begin{figure}[h]
 \caption{}\label{fig:swan2a}
$$\includegraphics[width=120mm]{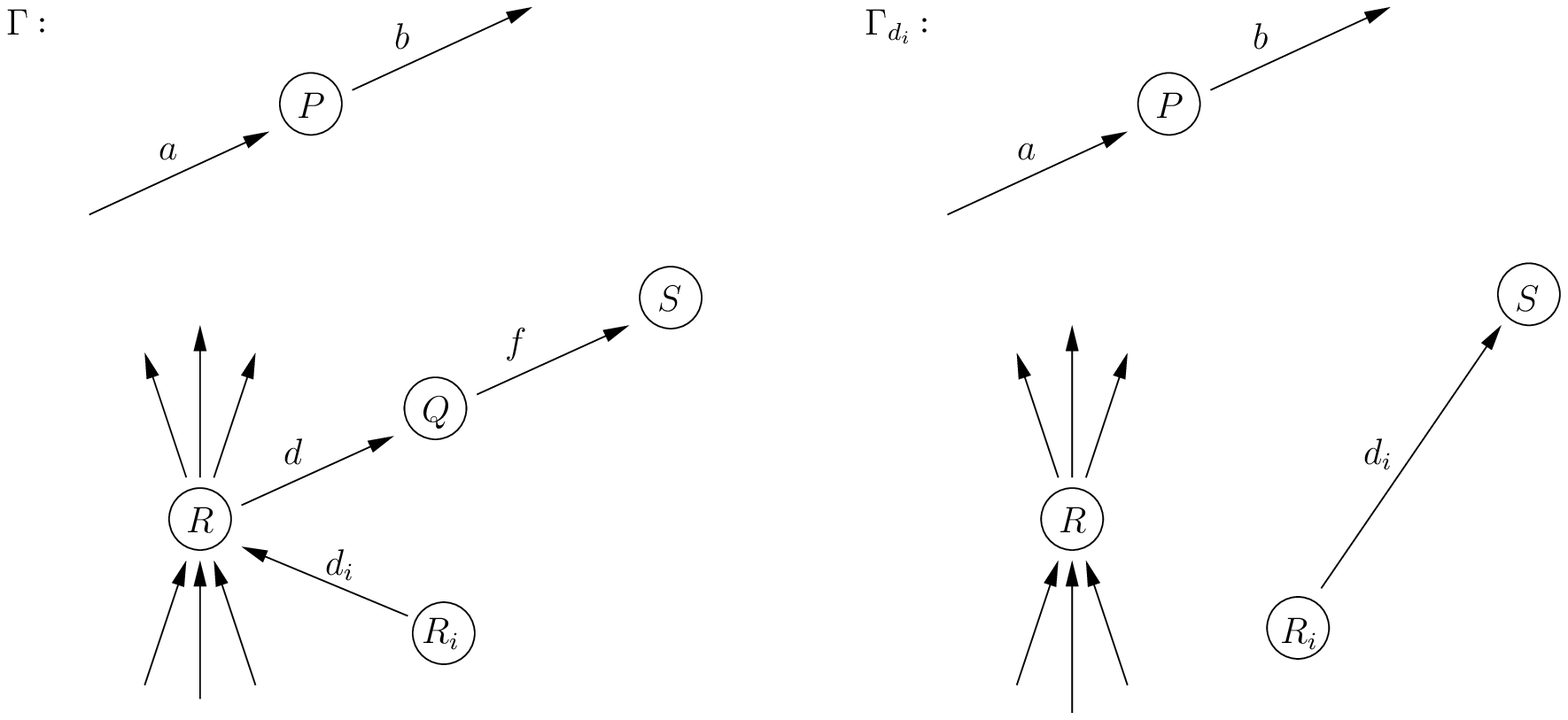}$$

\end{figure}

Case 2)  We are not in Case 1 and $P$ and $Q$ are not adjacent.    This is the same idea as in  Swan's case 2, page 371, combined with his induction argument of his case 1, page 370.   In the graph $\Gamma$  (See Figure \ref{fig:swan2a}),  we let $d_i$ be any of the edges that end at $R$ and let $g_j$ be any edge with initial vertex $R$.   We are trying to count  Eulerian paths on $\Gamma$ that contain the edges $a,b$ (in that order).   We partition those paths into two types:  Type (1)  consists of those paths that contain the sequence of edges $d_i,d,f$ or its reverse  $f^*,d^*,d_i^*$.   Type (2) consists of those paths that contain the sequence $g_j^*,d,f$ or its reverse $f^*,d^*,g_j$.  It is possible that there are no paths of type (2).    We will show that for each of these types the number of paths is even with half even path transformations of $\alpha$ and half odd.   That will prove the theorem in this case.

\begin{figure}[h]
 \caption{}\label{fig:swan2b}
$$\includegraphics[width=120mm]{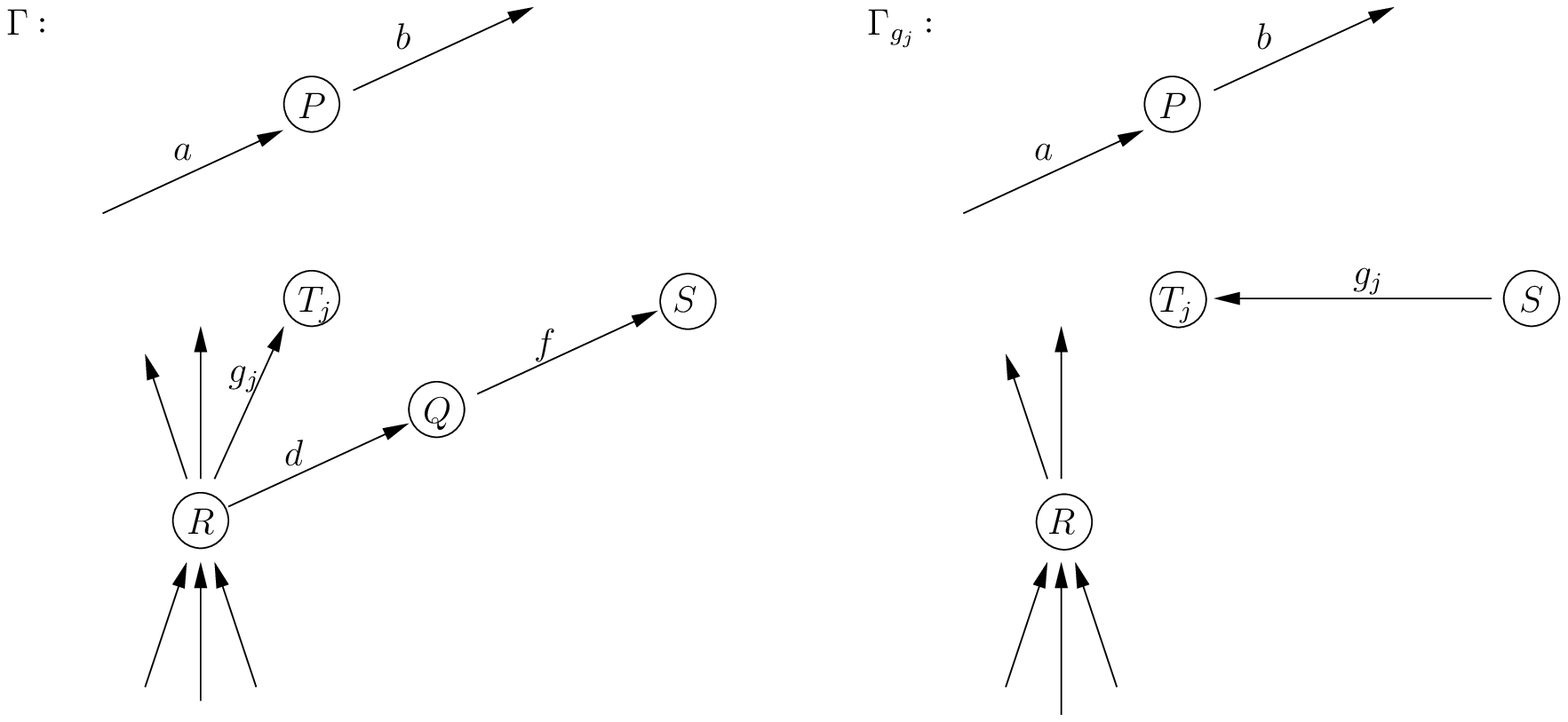}$$

\end{figure}

We first consider type (1).   For each such $d_i$ we create the new graph $\Gamma_{d_i}$  and new path $\alpha_{d_i}$ obtained by omitting the vertex $Q$ and the edges $d$ and $f$   and moving the endpoint of $d_i$ from $R$ to $S$ (See Figure \ref{fig:swan2a}).  Then the paths of type (1)   are in one to one correspondence with the paths on $\Gamma_{d_i}$ that have the same direction at $P$ as $\alpha_{d_i}$.   But  $\Gamma_{d_i}$  has $k-1$ vertices and $2k-4$ edges, so we can apply the induction hypothesis and infer that there are an even number of paths of type (1) , half even path transformations of $\alpha$, half odd.   Next we consider type (2).   We may assume there is such a path, say $\beta$.    Then for each  $g_j$ we create the new graph $\Gamma_{g_j}$  and new path $\alpha_{g_j}$ obtained by omitting the vertex $Q$ and the edges $d$ and $f$   and moving the initial point of $g_j$ from $R$ to $S$ (See Figure \ref{fig:swan2b}).  Then the paths on $\Gamma$ of type (2) are in one to one correspondence with the paths on $\Gamma_{g_j}$ that have the same direction at $P$ as $\beta$, that is the same direction as $\alpha$.    We then use the induction hypothesis on each $\Gamma_{g_j}$.    This finishes this case.

\begin{figure}[h]
  \caption{}\label{fig:swan3a}
$$\includegraphics[width=130mm]{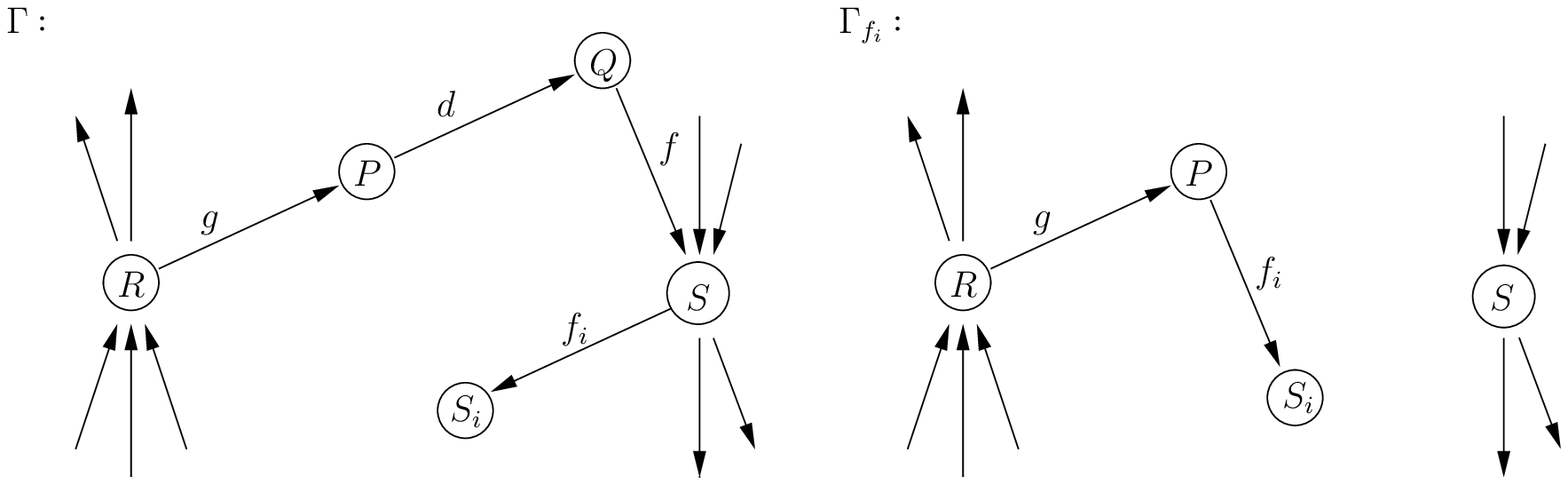}$$

\end{figure}

Case 3)  We are not in Case 1 and $P$ and $Q$ are adjacent.   In the graph $\Gamma$ (See Figure \ref{fig:swan3a}),  let $f_i$ be any of the edges that begin at $S$ and let $h_j$ be any edge that ends at $S$ (other than $f$).   We are trying to count Eulerian paths on $\Gamma$ that contain the edges $g,d$ (in that order),  which is the same as containing $g,d,f$.    Again there are two types:  Type (1) consists of those paths that contain $g,d,f,f_i$ and type (2) consists of those paths that contain $g,d,f,h_i^*$.  There may be no paths of type (2).    For type (1) we construct a new graph $\Gamma_{f_i}$ and new path $\alpha_{f_i}$ obtained by omitting the vertex $Q$ and the edges $d$ and $f$   and moving the initial vertex  of $f_i$ from $S$ to $P$ (See Figure \ref{fig:swan3a}).  Then the paths of type (1)   are in one to one correspondence with the paths on $\Gamma_{f_i}$ that have the same direction at $P$ as $\alpha_{f_i}$.   But  $\Gamma_{f_i}$  has $k-1$ vertices and $2k-4$ edges, so we can apply the induction hypothesis and infer that there are an even number of paths of type (1) , half even path transformations of $\alpha$, half odd.

\begin{figure}[h]
  \caption{}\label{fig:swan3b}
$$\includegraphics[width=130mm]{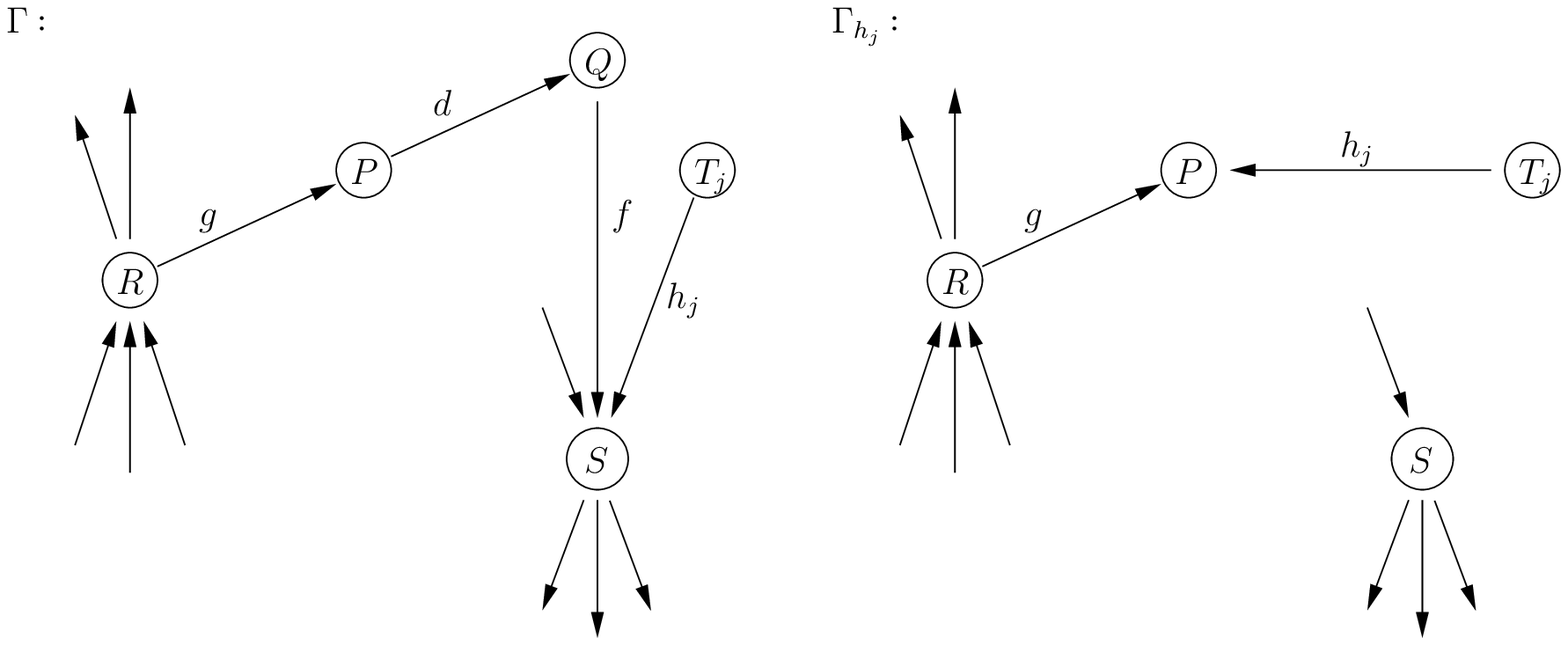}$$

\end{figure}

Next we consider type (2).   We may assume there is such a path, say $\beta$.    Then for each  $h_j$ we create the new graph $\Gamma_{h_j}$  and new path $\alpha_{h_j}$ obtained by omitting the vertex $Q$ and the edges $d$ and $f$   and moving the end point of $h_j$ from $S$ to $P$ (See Figure \ref{fig:swan3b}).  Then the paths on $\Gamma$ of type (2) are in one to one correspondence with the paths on $\Gamma_{h_j}$ that have the same direction at $P$ as $\beta$, that is the same direction as $\alpha$.    We then use the induction hypothesis on each $\Gamma_{h_j}$.    This finishes this case.\qed

Proof of Theorem \ref{swanlike}:   The argument of statement (3) of page 369 applies to our case and thus shows that we may assume the number of edges is exactly $2k-2$.   Let $\alpha$ be a path from $P$ to $Q$.  If $P=Q$ then by Proposition \ref{two vertices}  there is a vertex $R$ of order two and by Proposition \ref{loops}  there are an even number of path transformations, half odd and half even, that have the same direction as $\alpha$ at $R$.   But every path transformation is either one of these or the reverse of one of these, so the result holds in this case.

If $P\not= Q$, then we construct a new graph $\Gamma^\prime$  by adding one vertex $R$ and an edge $a$ from $Q$ to $R$ and an edge $b$ from $R$ to $P$.   Let $\alpha^\prime$ be the loop on $\Gamma^\prime$ that starts at $R$, then goes to $P$ then uses $\alpha$ and then goes back to $R$.  Clearly the paths from $P$ to $Q$ are in one to one correspondence with the paths  on $\Gamma^\prime$ that have the same direction at $R$ as $\alpha^\prime$.   We are therefore done by Proposition \ref{loops}.\qed

\section{Asymptotics of the codimension sequence}

In this section we analyze the asymptotic behavior of the codimension growth of the graded identities of $M_k(\CC)$ equipped with a crossed product grading and the transpose involution.
Our analysis follows the exposition of \cite[Section 3]{HN} and uses some techniques of \cite{RS}. Namely, we first establish a 1-1 correspondence between the equivalence classes of *-graded monomials and directed labeled graphs with an undirected Eulerian path, and thus measure the codimension by the number of these graphs. To obtain the asymptotics for the number of these graphs we count of the number of graphs of a more general kind and show the two asymptotics are the same.

To begin let $G = \{e=g_1, g_2, \ldots, g_k\}$ be a group of order $k$.
Let ${\bf g} = {\bf g}(k,n)$ denote a directed graph on $k$ vertices $\{g_1, g_2, \ldots, g_k\}$ labeled by the elements of the group and with $n$ edges labeled by the positive integers $\{1,2,\dots n\}$.

Recall that an Eulerian path in an undirected graph is a path that uses each edge exactly once. We call an Eulerian path a cycle if it starts and ends at the same vertex. We say that a directed graph has an Eulerian path or cycle if its underlying undirected graph does so.
We let $M_k(n)$ denote the set of all directed graphs ${\bf g}(k,n)$ that have an Eulerian path starting at the vertex $e$ and let $m_k(n)$ denote the order of $M_k(n)$.

 Recall from the Introduction that the $G$--graded $n$-codimension $c^G_n(A)$ of the algebra $A=M_k(\CC)$ is the dimension of the space $\displaystyle \frac{P^G_n}{P^G_n\cap I(G,*)}$, where $P^G_n$ is the $\qq-$vector space spanned by the monomials
$x^{\epsilon_n}_{i_1,g_1}x^{\epsilon_n}_{i_2,g_2}\dots x^{\epsilon_n}_{i_n,g_n}$ in the free algebra $\qq\{x_{i,g},x^*_{i,g}|i\geq 1, g\in G\}$, and $I(G,*)$ is the $(T,*)$--ideal of $G$--graded $*$--identities.

% As in \cite{HN}, there is
%one-to-one correspondence between the equivalence classes of monomials in $\displaystyle \frac{P^G_n}{P^G_n\cap I(G,*)}$ and the directed labeled graphs ${\bf g}(k,n) \in M_k(n)$.

\begin{thm}\label{independent}
Let $G$ be a group of order $k \geq 2$.
Then the graded $n$-codimension $c^G_n(A) = m_k(n)$. In particular, $c^G_n(A)$ does not depend on the group $G$.
\end{thm}

\proof By Corollary \ref{cor:binomials} the space $\displaystyle \frac{P^G_n}{P^G_n\cap I(G,*)}$ is generated by the the equivalence classes of strongly multilinear monomials.   We have seen that any such monomial  gives rise to a graph with an Eulerian path starting at the vertex $e$ and by Theorem \ref{thm:the same graph} two monomials are in the same class if and only if they define the same graph.   Conversely any Eulerian path starting at the the vertex $e$ in a graph ${\bf g}(k,n)$ in $M_k(n)$ clearly gives rise to a strongly multilinear monomial.  \qed

Because the $G$--graded $n$-codimension $c^G_n(A)$ depends on the order of the group $G$ only, we denote it $c_k(n)$.

We now describe the more general graphs we will need.

Let ${\bf g}(k,n)$ be a  directed, labeled graph with vertices the elements of $G$.  We define the {\it degree} of a vertex $g_i$ to be the number of edges having $g_i$ as their initial or terminal point. An edge may be a loop from $g_i$ to itself, and then it contributes twice to the degree of $g_i$.

It is easy to characterize an Eulerian path in such a graph. Namely, the graph $\bf g={\bf g}(k,n)$ has an Eulerian path from $g_i$ to $g_j$ if and only if

1. $\bf g$ is either connected, or the union of a connected subgraph and isolated points $g_l$ of $\deg(g_l) = 0$, where $l \neq i,j$.

2. for all $l \neq i,j$,  $\deg(g_l)$ is even,

3. if $i=j$, then $\deg(g_i)$ is even,

4. if $i \neq j$, then $\deg(g_i)$ and $\deg(g_j)$ are odd.

We say that a graph is {\it weakly connected} if it satisfies condition 1. We say that a graph is {\it strongly disconnected} if it is not weakly connected, i.e. has at least two non-trivial connected components. We say that a graph $\bf g$ has an Eulerian {\it pseudo-path} from $g_i$ to $g_j$ if it satisfies conditions 2-4, but is not necessarily weakly connected.  We say $\bf g$ is balanced  if   it has an Eulerian pseudo-cycle but is not necessarily weakly connected.

Let $\Gamma_k(n)$ be the set of all labeled, directed graphs ${\bf g}(k,n)$ on the vertices of $G$ which have an Eulerian pseudo-path from the vertex $e$ to the vertex $g_i$, for some $1\leq i \leq k$, and denote $|\Gamma_k(n)| = \gamma_k(n)$. Let $B_k(n)$ be the set of all balanced, directed graphs ${\bf g}(k,n)$, and denote $|B_k(n)| = b_k(n)$.

There is a simple formula relating  the number of graphs in the sets $\Gamma_k(n)$ and $B_k(n)$:

\begin{lem}\label{balanced-unbalanced}
$$\gamma_k(n) = \frac{1}{k}b_k(n+1).$$
\end{lem}

\proof Let $\widetilde{P}^{g_i}_k(n+1)$ be the set of all balanced graphs ${\bf g}(k,n+1)$ such that the endpoint of the edge labeled by $n+1$ is $g_i$. Clearly, there is 1-1 correspondence between any two sets $\widetilde{P}^{g_i}_k(n+1)$ and $\widetilde{P}^{g_j}_k(n+1)$ induced by the multiplication the node labels by $g_i^{-1}g_j$. Thus all these sets are of the same size and hence we have \[\displaystyle |\widetilde{P}^e_k(n+1)| = \frac{1}{k}|P_k(n+1)|.\]

Also, there is a 1-1 correspondence between the sets ${\widetilde{P}^e_k(n+1)}$ and $\Gamma_k(n)$: any graph ${\bf g}(k,n)$ having an Eulerian pseudo-path from the vertex $e$ to the vertex $g_i$, for some  $1\leq i \leq k$, can be completed in a unique way to a balanced graph ${\bf g}(k,n+1) \in \widetilde{P}^e_k(n+1)$ by adding the edge $n+1$ from $g_i$ to $e$. Thus we have
\[|\Gamma_k(n)| = |\widetilde{P}_k(n+1)| = \frac{1}{k}|P_k(n+1)|,\]
and the lemma follows.\qed

The next step is to count the number of balanced graphs on a fixed number of vertices $k$ with an arbitrary number of edges $n$.

\begin{lem}
The number of balanced labeled graphs ${\bf g}(k,n)$ is given by
\begin{equation}\label{multinomial}
\displaystyle b_k(n) = \sum_{n_1 + n_2 + \ldots + n_k = n}  \binom{2n}{2n_1\ 2n_2\ \dots \ 2n_k}.
\end{equation}
\end{lem}

\proof Given a set of $n$ edges labeled by the numbers $\{1,2, \ldots, n\}$ consider the set of heads and tails of all the edges. This is a set of cardinality $2n$. Let's call it the set of {\it the ends} of the edges. One constructs a unique labeled graph by choosing a vertex to which each one of the $2n$ ends is connected. The resulting graph is balanced if the number of ends at each vertex is even. Thus, for any partition $n_1 + n_2 + \ldots + n_k = n$ one chooses $2n_1, 2n_2, \dots, 2n_k$ ends connected to the vertices $g_1, g_2, \dots, g_k$  respectively. This gives the equation (\ref{multinomial}). \qed

The number $b_k(n)$ has the following properties:
\begin{lem} \label{recursive}
\begin{equation}\label{recursive 1}
b_k(n) = \sum_{i=0}^{n} \binom{2n}{2i} b_{k-1}(i).
\end{equation}
\begin{equation}\label{recursive 2}
b_{k_1+k_2}(n) = \sum_{i=0}^{n} \binom{2n}{2i} b_{k_1}(i)b_{k_2}(n-i).
\end{equation}
\end{lem}
\proof To prove equation (\ref{recursive 1}), choose $2n_1$ ends that are connected to the vertex $e$. This can be done in $\displaystyle \binom{2n}{2n_1}$ ways. Now connect the remaining $2n-2n_1$ ends to the vertices $g_2, \dots, g_k$ in $b_{k-1}(n-n_1)$ ways. Thus
$$b_k(n) = \sum_{n_1=0}^{n} \binom{2n}{2n_1} b_{k-1}(n-n_1) = \sum_{n_1=0}^{n} \binom{2n}{2n-2n_1} b_{k-1}(n-n_1) = \sum_{i=0}^{n} \binom{2n}{2i} b_{k-1}(i).$$

Now, let $k = k_1 + k_2$. Choose $2i$ ends that are connected to the vertices $g_1, \dots, g_{k_1}$. This can be done in $\displaystyle \binom{2n}{2i}$ ways. Once the ends are chosen, there are $b_{k_1}(i)$ ways to connect them the $k_1$ vertices. The remaining $2n-2i$ ends can be connected to the remaining $k_2$ vertices in $b_{k_2}(n-i)$ ways. This proves equation (\ref{recursive 2}). \qed

The next theorem gives the asymptotic expression for the number of balanced graphs.

\begin{thm} \label{Richmond} Let $k$ be an integer $\geq 2$. Then, as $n \rightarrow \infty$, we have

%\begin{equation}\label{multinomial asymptotics} (a) \hskip 1 in

%end{equation}

$$\displaystyle b_k(n) = \sum_{n_1 + n_2 + \ldots + n_k = n}  \binom{2n}{2n_1\ 2n_2\ \dots \ 2n_k} \sim  \frac{1}{2^{k-1}} k^{2n}.$$

\end{thm}

\proof   We follow the computations in \cite[Section 3]{RS} and first consider a single  multinomial

$$\binom{2n}{2n_1\ 2n_2\ \dots \ 2n_k}$$

where $n_1 + n_2 + \cdots + n_k = n$.

%$$\displaystyle \binom{2n}{2n_1\ 2n_2\ \dots \ 2n_k} \sim k^{2n}\exp\left(-k\sum_{1\leq i \leq k}x_i^2\right)(4\pi n)^{\frac{1-k}{2}}k^{\frac{k}{2}}.$$

We define $x_i$ for ${1\leq i \leq k}$  by $\displaystyle n_i = \frac{n}{k} + x_i\sqrt{n}$.    It follows that  $\sum_{1\leq i \leq k} x_i = 0$.

We first estimate $\prod_{1\leq i \leq k} (2n_i)!$ using  Stirling's formula written in the following form:

$$n! = e^{n \log n - n}\sqrt{2\pi n}(1 + O(n^{-1})),$$

and Taylor's formula:

$$\log (1+y) = y -\frac{y^2}{2} + O(y^3).$$

We obtain
\begin{equation*}
\begin{array}{lllll}
  \log 2n_i & = & \displaystyle \log (\frac{2n}{k} + 2x_i\sqrt{n})
   & = & \displaystyle \log \left(\frac{2n}{k}\left(1 + \frac{x_i k}{\sqrt{n}}\right)\right) \\
 & & & = & \displaystyle \log \frac{2n}{k} + \log \left(1 + \frac{x_i k}{\sqrt{n}}\right) \\
 & & & = & \displaystyle \log \frac{2n}{k}  + \frac{x_i k}{\sqrt{n}} - \frac{x^2_i k^2}{2n} + O(x^3_i n^{-3/2}).
\end{array}
\end{equation*}

\begin{equation*}
\begin{array}{lll}
 2n_i \log 2n_i & = & \displaystyle \left(\frac{2n}{k} + 2x_i\sqrt{n}\right)\left(\log \frac{2n}{k}  + \frac{x_i k}{\sqrt{n}} - \frac{x^2_i k^2}{2n} + O(x^3_i n^{-3/2})\right)\\
   & = & \displaystyle \left(\frac{2n}{k} + 2x_i\sqrt{n}\right)\log \frac{2n}{k} + 2x_i\sqrt{n} + x_i^2k + O(x_i^3 n^{-1/2}).
\end{array}
\end{equation*}

\begin{equation*}
\begin{array}{lll}
 2n_i \log 2n_i - 2n_i & = & \displaystyle \left(\frac{2n}{k} + 2x_i\sqrt{n}\right)\log \frac{2n}{k} - \frac{2n}{k} + x_i^2k + O(x_i^3 n^{-1/2}).
\end{array}
\end{equation*}

\begin{equation*}
\begin{array}{lll}
\sum_{1\leq i \leq k} \left( 2n_i \log 2n_i - 2n_i \right) & = & \displaystyle 2n \log \frac{2n}{k} - 2n + k \sum_{1\leq i \leq k} x_i^2 + O(x_i^3 n^{-1/2}).
\end{array}
\end{equation*}

We now make the additional assumption that $|x_i|\leq n^\epsilon$ where $0< \epsilon < 1/6$.  We then obtain

\begin{equation*}
\begin{array}{lll}
\sum_{1\leq i \leq k} \left( 2n_i \log 2n_i - 2n_i \right) & = & \displaystyle 2n \log \frac{2n}{k} - 2n + k \sum_{1\leq i \leq k} x_i^2 + O(n^{-1/2 + 3\epsilon}).
\end{array}
\end{equation*}

Now
\begin{equation*}
\begin{array}{lll}
\prod_{1\leq i \leq k} (2n_i)! & = &  \displaystyle \prod_{1\leq i \leq k} \left(\frac{2n}{k} + 2x_i\sqrt{n}\right)!  \\
& = & \displaystyle \prod_{1\leq i \leq k}  \exp \left(2n_i \log 2n_i - 2n_i\right) \sqrt{4\pi n_i}(1 + O(n^{-1})) \\
& = & \displaystyle \exp \left( \sum (2n_i \log 2n_i - 2n_i) \right) \prod_{1\leq i \leq k} \sqrt{4\pi \left(\frac{n}{k} + x_i\sqrt{n}\right)} (1 + O(n^{-1})) \\
& = & \displaystyle \exp \left( 2n \log \frac{2n}{k} - 2n + k \sum x_i^2 + O(n^{-1/2 + 3\epsilon})\right) \left( \frac{4\pi n}{k}\right)^{k/2}(1 + O(n^{-1/2})),
\end{array}
\end{equation*}
and so for $\epsilon < 1/6$ we obtain the asymptotic formula
\begin{equation*}
\begin{array}{lll}
\prod_{1\leq i \leq k} (2n_i)! & \sim & \displaystyle \exp \left( 2n \log \frac{2n}{k} - 2n + k \sum_{1\leq i \leq k} x_i^2 \right) \left( \frac{4\pi n}{k}\right)^{k/2}.
\end{array}
\end{equation*}

Recall that
$$\displaystyle (2n)! \sim \exp \left( 2n \log 2n - 2n \right) \left(4\pi n\right)^{1/2},$$
and then
$$\displaystyle \binom{2n}{2n_1\ 2n_2\ \dots \ 2n_k} \sim \frac{(2n)!}{\prod (2n_i)!} = \exp \left( 2n \log k - k \sum_{1\leq i \leq k} x_i^2 \right) k^{k/2} \left(4\pi n\right)^{(1-k)/2}.$$

We now consider the sum $\displaystyle \sum_{n_1 + n_2 + \ldots + n_k = n}  \binom{2n}{2n_1\ 2n_2\ \dots \ 2n_k}$ which can be approximated by the multiple integral
\begin{equation}\label{integral 1}
\displaystyle k^{2n}(4\pi n)^{\frac{1-k}{2}}k^{\frac{k}{2}}\int_0^n \int_0^n \dots \int_0^n \exp\left(-k\sum_{1\leq i \leq k}x_i^2\right) dn_1 dn_2 \dots dn_{k-1}.
\end{equation}

Since  $\displaystyle dn_i = \sqrt{n} dx_i$ we get (cf. \cite[Equation 10]{RS}):
$$\displaystyle \int_0^n \int_0^n \dots \int_0^n \exp\left(-k\sum_{1\leq i \leq k}x_i^2\right) dn_1 dn_2 \dots dn_{k-1} \sim $$ \begin{equation}\label{integral 2}
n^{\frac{k-1}{2}}\int_{-\infty}^\infty \int_{-\infty}^\infty \dots \int_{-\infty}^\infty \exp\left(-k\sum_{1\leq i \leq k}x_i^2\right) dx_1 dx_2 \dots dx_{k-1}
\end{equation}

The integrand is asymptotic to our desired quantity only if $|x_i|\leq n^{\epsilon}$ but outside this region the integrand is exponentially small.
It is found in \cite[Equation 12]{RS} that
\begin{equation}\label{integral 3}\int_{-\infty}^\infty \int_{-\infty}^\infty \dots \int_{-\infty}^\infty \exp\left(-k\sum_{1\leq i \leq k}x_i^2\right) dx_1 dx_2 \dots dx_{k-1} = \pi^{\frac{k-1}{2}} k^{-\frac{k}{2}}.
\end{equation}
Combining the above we get:
$$\displaystyle b_k(n)=\sum_{n_1 + n_2 + \ldots + n_k = n} \binom{2n}{2n_1\ 2n_2\ \dots \ 2n_k} \sim k^{2n}(4\pi n)^{\frac{1-k}{2}}k^{\frac{k}{2}}n^{\frac{k-1}{2}} \pi^{\frac{k-1}{2}} k^{-\frac{k}{2}} =  k^{2n} 2^{1-k}.$$
This completes the proof of the theorem. \qed

Next we show that the numbers $m_k(n)$ and $\gamma_k(n)$ have the same asymptotics:

\begin{prop}\label{the same asymptotics}
$$ \frac{m_k(n)}{\gamma_k(n)} \rightarrow 1, \ \mbox{as} \ n \rightarrow \infty.$$
\end{prop}

\proof The proof is similar to the proof of \cite[Proposition 14]{HN}.

Let $sc_k(n)$, $sd_k(n)$, and $sd^j_k(n)$ denote the number of connected graphs, strongly disconnected graphs, and strongly disconnected graphs with a connected $e$-component on a given subset of the vertices of size $j$ in $\Gamma_k(n)$, respectively.
Then, as in \cite{HN}, $$ sd_k(n) = \sum_{j=1}^{k-1} \binom{k-1}{j-1} sd^j_k(n),$$
and $$sd^j_k(n) = \sum_{i=j-1}^{n-1} \binom{n}{i} sc_{j}(i) b_{k-j}(n-i),$$
and, using Lemma \ref{balanced-unbalanced}, we get
$$sd^j_k(n) \leq \sum_{l=1}^{n} \frac{l}{n+1} \binom{n+1}{l} \frac{1}{j} b_{j}(l) b_{k-j}(n+1-l).$$

%Here is the only place where the proof differs from the original one.

We claim that for any $n$ and $1\leq l\leq n$ holds $$\displaystyle l \binom{n+1}{l} < \binom{2n+2}{2l}.$$
Indeed, since $1\leq l\leq n$ we have
$$l \binom{n+1}{l} < (n+1) \binom{n+1}{l} \leq  \binom{n+1}{l}\binom{n+1}{l} \leq \binom{2n+2}{2l}.$$
The last inequality follows from the observation that one of the strategies to choose $2l$ out of $2n+2$ is to choose the first $l$ out of the first $n+1$, and, independently, the remaining $l$ out of the remaining $n+1$, and this strategy is not exhaustive.

It then follows by Lemma \ref{recursive} that
          $$sd^j_k(n) \leq  \frac{1}{n+1}\sum_{l=0}^{n+1} \binom{2n+2}{2l} b_{j}(l) b_{k-j}(n+1-l) = \frac{1}{n+1}b_{k}(n+1).$$
Hence, as $n \rightarrow \infty$, we have
          $$\frac{sd^j_k(n)}{\gamma_k(n)} \leq \frac{\frac{1}{n+1}b_{k}(n+1)}{\frac{1}{k}b_k(n+1)} = \frac{k}{n+1} \rightarrow 0,$$ and
          $$\displaystyle \frac{sd_k(n)}{\gamma_k(n)} = \sum_{j=1}^{k-1} \binom{k-1}{j-1} \frac{sd^j_k(n)}{\gamma_k(n)} \rightarrow 0,$$
and finally
$$\displaystyle \frac{m_k(n)}{\gamma_k(n)}  = \frac{\gamma_k(n) - sd_k(n)}{\gamma_k(n)}\rightarrow 1.$$
This completes the proof of the proposition. \qed

We are now ready to prove the main theorem of this section:

\begin{thm}\label{asymptotic} Let $G$ be a group of order $k \geq 2$. Then, as $n \rightarrow \infty$, the $G$--graded $n$-codimension
of $M_k(\CC)$ equipped with the $G$--crossed product grading and the transpose involution is
$$\displaystyle c_k(n) \sim \frac{k}{2^{k-1}} k^{2n}.$$
\end{thm}

\proof By Proposition \ref{independent} we have $c_k(n) = m_k(n)$. Hence, by Proposition \ref{the same asymptotics}, $c_k(n) \sim \gamma_k(n)$. By Lemma \ref{balanced-unbalanced}, $\displaystyle c_k(n) \sim \frac{1}{k}p_k(n+1)$.  Finally by Theorem \ref{Richmond},

$$\displaystyle c_k(n) \sim \frac{1}{k}\left(\frac{1}{2^{k-1}} k^{2(n+1)}\right) = \frac{k}{2^{k-1}} k^{2n}.$$
\qed

In case $k=2$ we give an explicit value of the $n$-codimension:

\begin{thm}\label{asymptotic-case-2} Let $C_2$ be a cyclic group of order $2$. Then the $C_2$--graded $n$-codimension of $M_2(\CC)$
equipped with the $C_2$--crossed product grading and the transpose involution is
        $$\displaystyle c_2(n) = 4^n - 2^n + 1.$$
In particular, as $n \rightarrow \infty$, we have
        $$\displaystyle c_2(n) \sim 2^{2n}.$$
\end{thm}

\proof Let $C_2 = \langle e, \sigma \rangle$ be a cyclic group of order 2.

Notice that any weakly connected directed graph on the vertices $e$ and $\sigma$ necessarily has an (undirected) Eulerian path. The total number of directed labeled graphs on two vertices with $n$ edges is $2^{2n}$ (this is the number of all possible distributions of $2n$ ends between two vertices). Among these graphs, exactly $2^n$ are disconnected. Among the disconnected graphs one (that having all the edges starting and ending at $e$) is weakly connected. Thus we have
$$\displaystyle c_2(n) = 2^{2n} - 2^n + 1.$$
The second statement of the theorem follows at once.\qed

\end{document}